\documentclass[11pt,a4paper]{article}

\usepackage[english]{babel}

\usepackage{amsmath}
\usepackage{amssymb}
\usepackage[dvips]{graphicx}
\usepackage{amsmath}
\usepackage{amssymb}
\usepackage{amsbsy}
\usepackage{amsfonts}
\usepackage{array}
\usepackage{verbatim}
 \usepackage{bm}
\usepackage{array}
\usepackage{booktabs}
\usepackage{caption}
\usepackage{tabularx}

\newtheorem{theorem}{Theorem}

\begin{document}

\title{Local Exponential Methods: a domain decomposition approach to
exponential time integration of PDEs}

\author{Luca Bonaventura}

\maketitle

\begin{center}
{\small
MOX -- Modelling and Scientific Computing \\
Dipartimento di Matematica, Politecnico di Milano \\
Via Bonardi 9, 20133 Milano, Italy\\
{\tt luca.bonaventura@polimi.it}\\
}
\end{center}

\date{}


\vspace*{0.5cm}


\vspace*{0.5cm}

\abstract{A local approach to the time integration of PDEs by exponential methods is
proposed, motivated by theoretical estimates by A.Iserles on the decay of off-diagonal
terms in the exponentials of sparse matrices. An overlapping domain decomposition technique
is outlined, that allows to replace the computation of a global exponential matrix
by a number of independent and easily parallelizable local problems. Advantages and potential
problems of the proposed technique are discussed. Numerical experiments on simple, yet relevant
model problems show that the resulting method allows to increase  computational efficiency 
with respect to standard implementations of exponential methods.}

\pagebreak


\section{Introduction}
\label{intro} 
\indent

\begin{figure}[htbc]
\begin{center}
\includegraphics[height=0.40\textheight,width=0.95\textwidth]{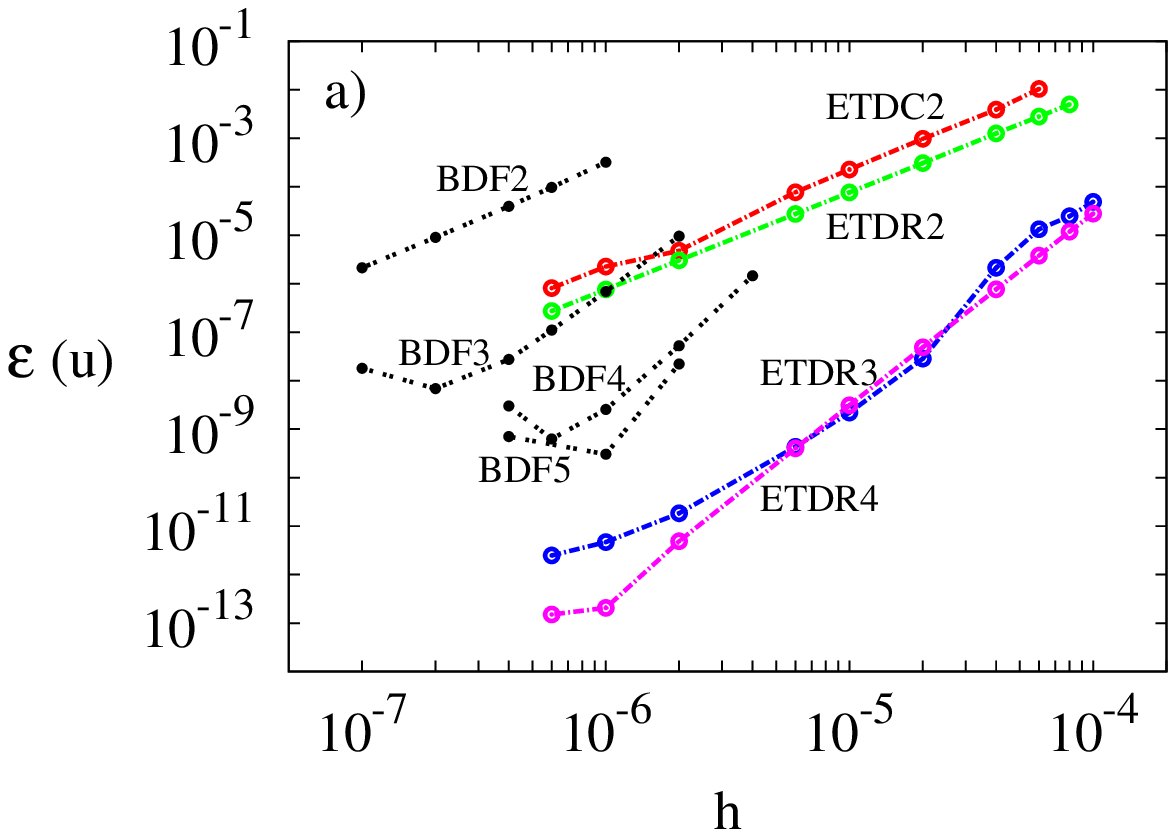}(a)
\includegraphics[height=0.40\textheight,width=0.95\textwidth]{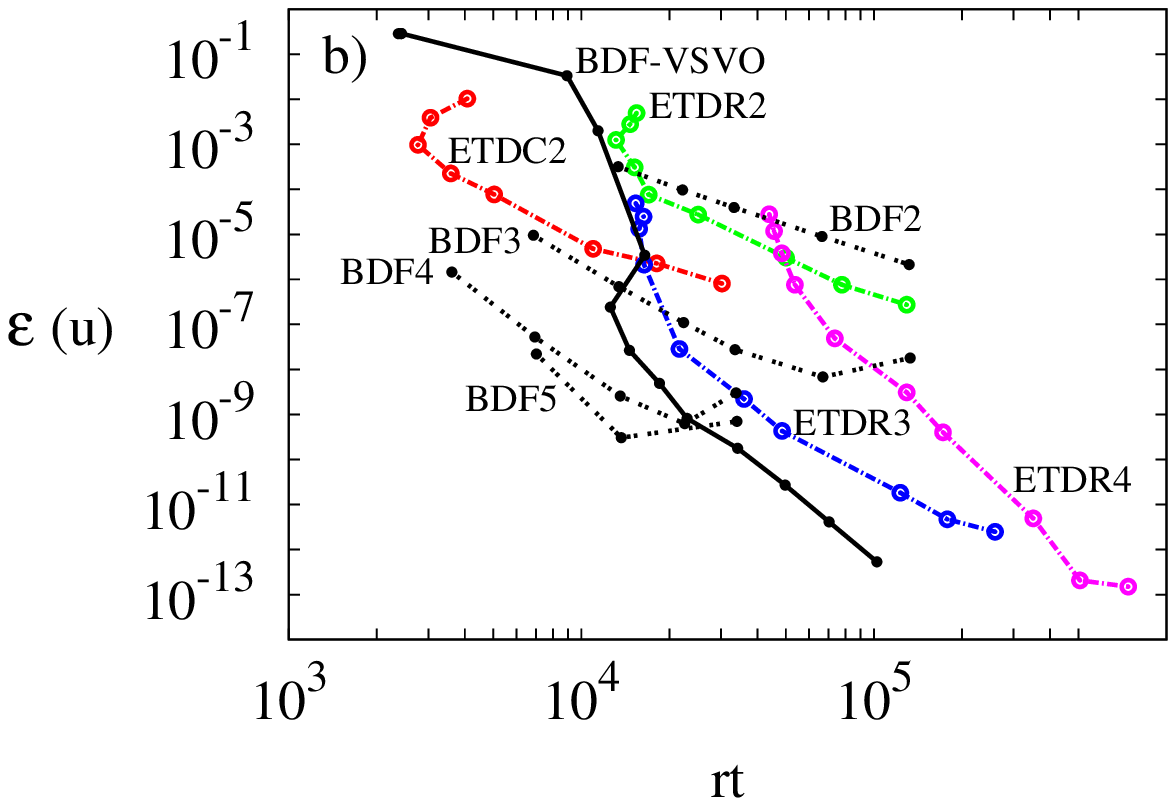}(b)
\end{center}
\caption{Errors  as a function of time step (a) and computational cost (b) for various IMEX and exponential methods. Reproduced from
\cite{garcia:2014} with the consent of the authors.}
\label{imex}
\end{figure}

The application of exponential time integration methods (EM) to the time discretization of partial differential equations
has been the focus of increasing attention  over the last two decades. A recent review of these methods
is provided e.g. in \cite{hochbruck:2010}. EM   allow  to eliminate almost entirely
the time discretization error in the linear case and to reduce it substantially in many nonlinear cases. 
EM possess stability properties that make them competitive with standard
stiff solvers. Furthermore,  when wave propagation problems are considered,
 they also allow to represent faithfully  even the fastest linear waves,
which are usually damped and distorted by conventional implicit and semi-implicit techniques.
After the seminal papers \cite{gallopoulos:1992}, \cite{saad:1992}, 
various kinds of exponential  methods have been employed
by a  number of authors in conjunction with different
time discretization techniques, see e.g. \cite{archibald:2011}, \cite{beylkin:1998}, \cite{garcia:2014}, \cite{madaule:2014},
\cite{martinez:2009}, \cite{schulze:2009}.
 In \cite{garcia:2014}, an extensive comparison
between IMEX and exponential methods has been carried out. Both kinds of time discretizations have been coupled to  
 a spectral spatial discretization meant for applications to mantle convection problems, but also very similar to
 numerical techniques widely used in numerical weather prediction and climate modelling.
One of the conclusions of this work is highlighted by the graphs in figure \ref{imex}, that
are representative of a large number of numerical simulations carried out in that paper.

While superior in terms of the accuracy achieved when using  large time steps, exponential
methods appear much more expensive per time step than IMEX methods. Numerical
experiments carried out by the author with  semi-implicit methids employed in atmospheric
modelling point to the same  conclusions. The high computational cost is due to the effort in computing
the matrix exponential  by the Krylov space techniques of \cite{saad:1992}. Numerical experiments
carried out with other approaches for the exponential matrix computation,
 such as the Leja points interpolation employed in \cite{martinez:2009},
do not seem to improve things substantially in this respect.

These results motivate the quest for a more economical approach to the computation of the exponentials
of sparse matrices typically arising in the space discretization of PDEs.
This sparsity  corresponds to an effectively finite domain
of dependence   for the solution of both hyperbolic and parabolic problems. 
This locality in space of the exact solution of the PDE makes it appear illogical that the exact computation 
of the spatially discretized solution should be a global  operation,     involving all the degrees of freedom
of the problem, rather than only those that are close  to each other in physical space.
Such a heuristic consideration has a rigorous counterpart
in the results proven  in \cite{iserles:2001} by A.Iserles on the decay of off-diagonal
terms in the exponentials of sparse matrices. The bounds proven in \cite{iserles:2001} imply that,  given
a sparse matrix ${\bf A} $ representing the spatial discretization of a local, linear differential operator  and a 
vector ${\bf x} $ representing the discrete degrees of freedom, for any node $i $
the contribution to $ ({\exp({\bf A}\Delta t)} {\bf x})_i $ of nodes $j$ that are  sufficiently 
far from $i$ in physical space is small. How far and how small will obviously depend on the time step
 $ \Delta t  $ and on the speed of propagation of information in the PDE discretized by ${\bf A}.$ 
 Extensions of this result have been proposed more recently in \cite{benzi:2014}, \cite{benzi:2015b}, \cite{benzi:2015a}.

Motivated by these considerations, an overlapping domain decomposition technique
is proposed here,  that allows to replace the computation of a global exponential matrix
by a number of independent local matrices. The obvious advantage of such a Local Exponential Method (LEM)
 is that each local problem can be solved independently in parallel, thus  increasing the scalability
of the resulting time discretization technique. Furthermore,  if the  number
of degrees of freedom associated to each local domain is small enough, the local exponential matrices can be computed
by   Pad\'e approximation combined with scaling and squaring 
(see e.g. \cite{moler:2003} for a review of numerical methods for the computation of matrix exponentials)
 and can be actually stored in memory, thus bypassing the
problems that result from having to compute the action $ ({\exp({\bf A}\Delta t)} {\bf x}) $ rather than
the exponential matrix itself. The main drawback of the proposed approach is that, in each of the
local problems, some of the nodes also associated to other subdomains are updated just for the
sake of providing a buffer zone that makes the domain under consideration only marginally affected
by the far field. For hyperbolic problems, the size of this buffer zone depends on the Courant
number, which clearly implies that the method becomes increasingly inefficient in the limit
of very large Courant numbers. However, some situations exist in which the proposed approach 
seems competitive in spite of this limitation. 
In some linear problems relevant for many applications, such as the Schr\"odinger equation,
being able to store 
 the local matrices required for approximation of
the exponential matrix can significantly reduce the computational cost.
In high order finite element methods, 
high Courant numbers arise easily due to the large number of degrees of freedom per element
(see e.g. the related discussion and proposals in \cite{tumolo:2013}),
so that a technique that enables to run at Courant numbers of the order of the polynomial
degree employed with a more local approach could turn out to be useful. Furthermore, in many environmental
applications, such as numerical weather prediction and ocean modelling, strongly anisotropic
meshes are employed, with a vertical resolution that is often  two or three orders of magnitude
smaller than the horizontal one. This results in high vertical Courant numbers, that are often
addressed by directional splitting methods. The present approach would allow to achieve
the same goal by employing a horizontal domain decomposition approach with minimal
overlap among subdomains, such as
almost universally used for parallelization of this kind of models, while at the same time
avoiding \textit{ad hoc} solutions that rely on splitting and providing an efficient 
and robust way to solve  the corresponding fluid dynamics equations.

In section \ref{sparse}, the key results of \cite{iserles:2001} are reviewed and
applied to the spatial discretization of a simple model problem. In section \ref{dd}, an overlapping
 domain decomposition approach to exponential time integration methods is outlined
 and the potential advantages and  problems of the proposed technique are discussed.
 In section \ref{tests}, some numerical results on simple one and two dimensional
 problems are reported, that show that the proposed approach is able to attain the same
 accuracy level as either standard time discretization methods or global implementations
 of exponential integrators. 
 Some conclusions and perspectives for future work are outlined in section \ref{conclu}.

\section{Exponentials of sparse matrices and PDEs}
\label{sparse} 
\indent
The key theoretical result for the development proposed in this paper has been  proven in \cite{iserles:2001}
and will be briefly summarized here. 

\begin{theorem} \label{basic_sparse}
Consider a sparse, $s-$banded matrix ${\bf A}=(a_{i,j})  $
and assume that its non zero entries are bounded by
 $\max_{i,j}|a_{i,j}| \leq \rho.$  Denoting then
$ \exp({\bf A})=(e_{i,j}) $ one has 

\begin{eqnarray}  
 |e_{i,j}| &\leq & \Big (\frac{\rho s}{|i-j|} \Big )^{\frac {|i-j|}s}\Big [ {\rm e}^{\frac {|i-j|}s} - \sum_{k=0}^{|i-j|-1 }  \frac{(| i-j/s |)^k}{k!}\Big ]
\nonumber\\
  &\approx &\Big (\frac{\rho s}{|i-j|} \Big )^{\frac {|i-j|}s} \frac{(|i-j|/s)^{|i-j|}}{|i-j|!}  
  \end{eqnarray}
  \end{theorem}
  Some remarks on this result are necessary before proceeding to the application to spatial discretizations of PDEs.
  Firstly, in many exponential matrix it is necessary or convenient to use, rather than the exponential itself, the so
  called $\phi - $ functions, defined recursively as
  $$\phi_k(z) = \frac{\exp{(z)}-\sum_{l=0}^kz^l/l! }{z^k}.$$
  The function $\phi_1$ will be simply denoted as $\phi$ in the following. Although theorem \ref{basic_sparse} as stated in \cite{iserles:2001}
  refers strictly speaking to the exponential function only, the bounds given in the same paper
  on the size of the entries of powers ${\bf A}^m$ could be employed to derive similar decay estimates for the off diagonal
  terms of $\phi_k({\bf A}).$  Here, we do not pursue the rigorous extension of \ref{basic_sparse} to exponential
  related functions. Furthermore, theorem \ref{basic_sparse} is proven strictly speaking only for matrices with limited bandwith,
  while in spatial discretizations of multidimensional PDEs more complex sparsity patterns easily arise. Here, it will be again
  assumed heuristically that analogous estimates can be provided also in these more relevant cases, although it is
  clear that deriving general proofs may be far from easy.  Numerical results reported  in section \ref{tests} seem to provide heuristic evidence that 
   the required generalizations of \ref{basic_sparse} hold. 
  
  Given this \textit{caveat}, it is possible to argue that
theorem \ref{basic_sparse} has important implications for the application of exponential methods to the time discretization of PDEs.
We will consider as a model problem the linear advection diffusion equation and its associated initial and boundary value problem on the spatial domain ${\Omega} $
\begin{eqnarray}
&&\frac{\partial   c } {\partial t}=-{\bf a}(x) \cdot \nabla c  +\nu (x) \Delta   c , \ \ \ \ t \in[0,T]  \ \ \ x \in \Omega \nonumber \\
&&  c(x,0)=  c_0(x) \label{ibv}
\end{eqnarray}
with appropriate boundary conditions (taken here for simplicity to be time independent) defined by the linear operator
 ${\cal B} c=   g $ at $\partial \Omega.$ 
Here ${\bf a}(x) $ denotes the velocity field and $\nu (x) $ the non negative diffusivity.
 Denote  then by ${\cal M}$   a computational
mesh with minimum element size $\Delta x $ for the approximation of \eqref{ibv} by some consistent space discretization technique. Denote
by ${\bf u} $ the vector of the  associated discrete degrees of freedom and by ${\bf A} $ the matrix representing
the spatial discretization. 
Problem \eqref{ibv} will then be approximated by
\begin{eqnarray}
&&\frac{d {\bf u} } {d t} = {\bf A}{\bf u}+{\bf g} \nonumber \\
&&{\bf u} (0)={\bf u}_0, 
\label{spacedisc}
\end{eqnarray}
where now by a slight abuse of notation ${\bf g} $ denotes the non autonomous forcing resulting from
the spatial discretization of the boundary conditions ${\cal B}  c=  g$ and  ${\bf u}_0 $ denotes the approximation
of $ c_0(x) $ in the chosen finite dimensional setting.
The  solution of problem \ref{spacedisc} is given by
$$  {\bf u}(t) = \exp{({\bf A}t)}  {\bf u}_0 
+\int_0^t \exp{({\bf A}(t-s))} {\bf g} \ ds = {\bf u}_0 + t\phi\big (t{\bf A} \big ) \Big ({\bf A}{\bf u}+{\bf g}  \Big ).
$$
It is well known that this formula is the basis of all exponential methods, that in many cases
simply reduce to it for linear problems with constant forcing.
For any reasonable, consistent  spatial discretization, it will follow that
the entries of  $\Delta t{\bf A}$ will be bounded by
  $$ \rho = O \Big (C + \mu \Big ),  \ \ \ C=\max_{x  \in  \Omega} |{\bf a} (x)|\frac{\Delta t} {\Delta x}
  \ \ \ \  \mu=\max_{x  \in  \Omega} \nu(x) \frac{\Delta t} {\Delta x^2} , $$
  where $C,\mu$ denote, respectively, the maximum Courant number and the stability parameter of standard
  explicit discretizations of the heat equation.
  
\begin{figure}[htbc]
\begin{center}
\includegraphics[height=0.3\textheight,width=0.8\textwidth]{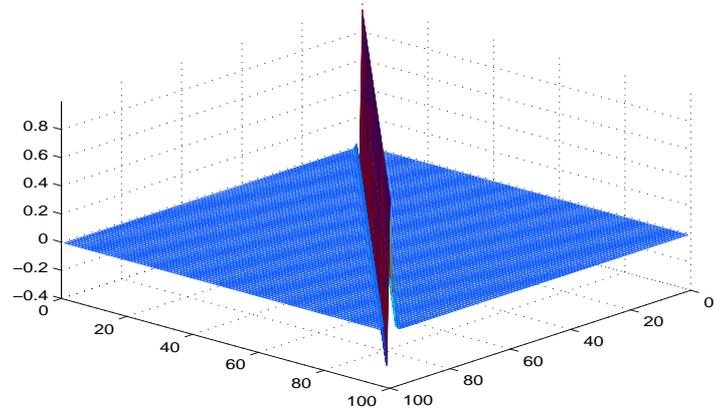}(a)
\includegraphics[height=0.3\textheight,width=0.8\textwidth]{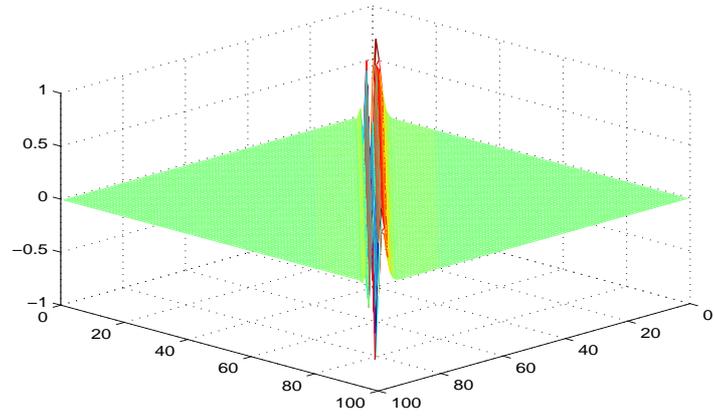}(b)
\includegraphics[height=0.3\textheight,width=0.8\textwidth]{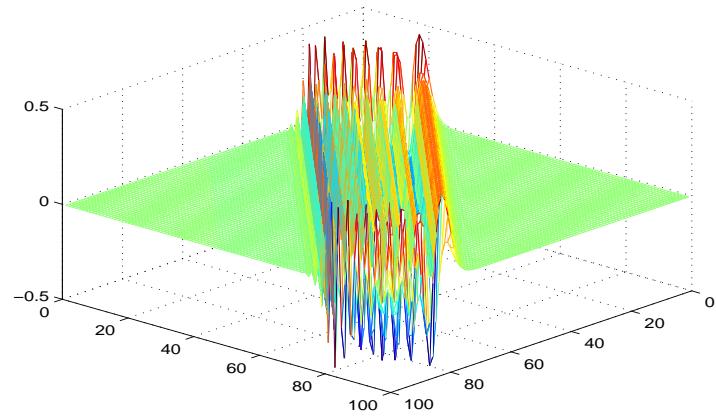}(c)
\end{center}
\caption{Visualization of    ${\exp}(\Delta t{\bf A})$ in the case of the centered finite difference approximation of 
advection in 1d at Courant numbers (a) 0.5, (b) 5, (c) 20.}
\label{adelta}
\end{figure}

  This implies that the degree of sparseness of  ${\exp}(\Delta t{\bf A})$ and of the related functions will depend on
  the magnitude of the usual stability parameters for explicit time discretizations.
  This can be visualized graphically by considering a simple centered finite difference
  approximation on a uniform mesh for the one dimensional case with $\mu=0$
  with homogeneous Dirichlet boundary conditions. The matrix  $\Delta t{\bf A}$
  is visualized in figure \ref{adelta}    for    Courant numbers $ 0.5,  5,  20, $ respectively.
  Entries away from the diagonal are non zero, but have an absolute value that decays rapidly
  for bands further away from the diagonal than the Courant number.
Even though the diffusion operator has an infinite domain of dependence in the continuous case, a similar pattern
is observed for simple discretizations of the diffusive terms. Strategies to obtain
optimal bounds for the entries of the exponential are proposed in \cite{iserles:2001} and the issue of how to estimate
rigorously an optimal bound in more general cases is definitely an open one. However, as shown
in the practical applications to  PDE problems  presented in section \ref{tests},
one can simply rely on the inspection of the standard stability parameters, like the Courant number,
to identify the size of the mesh region that is effectively contributing to the change in the solution
associated with a given node within one time step.
These considerations lead to the idea that the computation of the global exponential matrix required
by standard approaches to application of exponential integrators for PDEs can be replaced by the
computation of local matrices, associated to the exact solution of the restriction of \eqref{spacedisc}
to  appropriate subsets of the computational domain.  

\section{Local exponential methods: a domain decomposition approach to exponential time integration methods}
\label{dd} 
\begin{figure}[htbc]
\begin{center}
\includegraphics[height=0.5\textheight,width=0.7\textwidth]{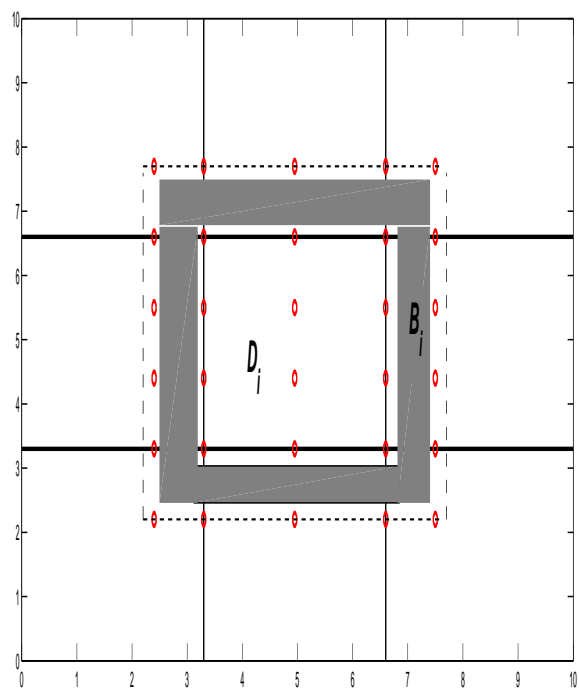}
\end{center}
\caption{Sketch of a domain ${\cal D}_i $ with the corresponding buffer region ${\cal B}_i $ highlighted in grey. }
\label{buffer}
\end{figure}

\indent
The theoretical results summarized in the previous section suggest a more
local approximation of exponential matrices for time discretization of linear PDEs.
One possible approach will be now be introduced and denoted shortly  as Local Exponential Method (LEM) in the following.
The first step consists in
decomposing the mesh in $D$ overlapping regions $${\cal M}= \bigcup_{i=1}^D{\cal M}_i  \ \ \   {\cal M}_i ={\cal D}_i \cup {\cal B}_i .$$
Here ${\cal D}_i $ denote non overlapping  domains, such that  ${\cal M}= \bigcup_{i=1}^D{\cal D}_i ,$
while ${\cal B}_i $  are boundary buffer zones surrounding each ${\cal D}_i. $ A visualization of one such region is   sketched in figure \ref{buffer}.

  Notice that, after the space discretization \eqref{spacedisc}, the mesh is only assumed to include
interior nodes, while the effect of boundary conditions is included in the forcing term.
One can then denote by ${\bf u}_{{\cal M}_i},$  ${\bf u}_{{\cal D}_i}$ the set of discrete degrees of freedom 
associated to ${\cal M}_i, {\cal D}_i, $ respectively and by ${\bf A}_{{\cal M}_i} $ the restriction of the matrix
${\bf A}$ to the nodes in ${\cal M}_i. $ 
Given these definitions, we now outline the LEM in the case of the simplest exponential method, i.e. the exponential
Euler method (see e.g. \cite{hochbruck:2010}). The extension to
any of the exponential methods described in the literature is immediate. 
Introduce a discrete set of time levels $ t_n,$ taken for simplicity to be equally spaced and such
that $t^{n+1}-t^n=\Delta t.$ Then,  for each $i=1,\dots,D,$  and $t^{n},$ 

\begin{itemize}
\item[1.]  consider the local problems restricted to $  {\cal M}_i $ 
\begin{equation}
\frac{d {\bf v}_{{\cal M}_i} } {d t} = {\bf A}_{{\cal M}_i}{\bf v}_{{\cal M}_i}+{\bf g}_{{\cal M}_i}(t),
\label{local_pb}
\end{equation}
where ${\bf v}_{{\cal M}_i}(t^n) $ is assumed to coincide  with ${\bf u}^n_{{\cal M}_i} $ 
and $ {\bf A}_{{\cal M}_i},{\bf g}_{{\cal M}_i}(t),$ are modified with respect to the simple restriction
of their global counterparts in order to impose Dirichlet boundary  conditions along the parts of the boundaries $\partial {\cal B}_i $
that belong to some domain $ {\cal D}_j $ with $j\neq i;$
\item[2.] compute
\begin{equation}
{\bf v}_{{\cal M}_i}^{n+1}= {\bf v}_{{\cal M}_i}^{n}
+\Delta t \phi(\Delta t {\bf A}_{{\cal M}_i}  ) \Big ({\bf A}_{ {\cal M}_i} {\bf v}^n_{ {\cal M}_i} +{\bf g}_{{\cal M}_i}(t^n) \Big );
 \label{lem}
 \end{equation}
\item[3.]   
define $ {\bf u}_{{\cal D}_i}^{n+1}= {\bf v}_{{\cal D}_i}^{n+1}, $ thus
 overwriting the degrees of freedom belonging to the buffer zones.
\end{itemize}
It is clear that, if the solution of each  local problem \eqref{local_pb} is to be a good approximation of the solution of
global problem \eqref{spacedisc}, the buffer regions must be chosen in such a way that the contribution 
to $  \phi(\Delta t {\bf A}_{{\cal M}_i}  )$ from nodes $k \notin {\cal M}_i $ is negligible.
As discussed in the previous section, for discretizations of the model problem \eqref{ibv} it should
be sufficient to choose a size that, given a value of $\Delta t, $ is related to the typical stability
parameters of explicit time discretizations. In the simple tests described in section \ref{tests}, on
cartesian meshes the size of the buffer regions is  taken to be (empirically) related to the
maximum between $C $ and $   \mu.$

  A number of obvious potential advantages of LEM are apparent. First of all, 
  the computation of the local exponential matrices is trivially parallel, thus leading to an algorithm
  that should scale much better on massively parallel machines than those requiring a 
  global communication step. Furthermore, the size of the exponential matrices to be computed will
  be smaller, thus implying that, if e.g. Krylov space techniques are employed for their computation,
 a smaller dimension of the Krylov space would be sufficient for their accurate approximation.
   Finally, for small enough subdomains ${\cal M}_i,$ the resulting local matrices  could be 
   computed   by simpler methods, such as direct Pad\'e approximation (see e.g. \cite{moler:2003}),
   and stored in the local memory. This contrasts with  the computation of the action of the exponential
   matrix that is necessary in standard implementations for large ODE systems deriving from spatial
   semidiscretization of PDEs. Storage of the local matrices would also allow to
  to increase the efficiency of methods such as the exponential
   Rosenbrock methods proposed in \cite{hochbruck:1997}, 
    by allowing to freeze the Jacobian matrix of the right hand side over a certain number of time steps.

   Some disadvantages of the proposed approach are also obvious. The degrees of freedom belonging 
   to the buffer regions only play an auxiliary role and would be updated at least twice
   (or more, if the corresponding nodes belong to more than two buffer regions). This implies that
   there is a computational overhead that is proportional to the size of such regions. Considering for simplicity the
   pure advection problem, this implies that the method is increasingly less efficient in the limit
   of increasing Courant number, which is exactly the limit in which exponential methods are most advantageous
   with respect to more standard ones. In which regimes the resulting algorithm could end up in being more efficient
   than standard ones is not obvious.
   However,  some situations can easily be identified in which the proposed approach 
seems competitive in spite of this limitation. In high order finite element methods, for example,
high Courant numbers arise easily due to the large number of degrees of freedom per element (see e.g. \cite{tumolo:2013}),
so that a technique that is able to run at Courant numbers of the order of the polynomial
degree employed with a more local approach could turn out to be useful. Furthermore, in many environmental
applications, such as numerical weather prediction and ocean modelling, strongly anisotropic
meshes are employed, with a vertical resolution that is often  two or three orders of magnitude
smaller than the horizontal one. This results in high vertical Courant numbers, that are often
addressed by directional splitting methods. The present approach would allow to achieve
the same goal by employing a horizontal domain decomposition approach with minimal
overlap among subdomains, such as
almost universally used for parallelization of this kind of models, while at the same time
avoiding \textit{ad hoc} solutions that rely on splitting and providing an efficient 
and robust way to solve  the corresponding fluid dynamics equations. 
The preliminary numerical results reported in section \ref{tests} support the view
that an acceptable trade off between locality and efficiency is feasible
and motivate further investigation of the application of LEM to fluid dynamics and wave propagation problems.
Furthermore,
if environmental models with complex physical parameterizations are considered,
LEM would provide  a completely local approach to include these terms while
maintaining high order accuracy in time without extra computational costs.
Indeed, if the second order exponential Rosenbrock method  is considered,
second order accuracy in time would be attainable with a single evaluation of the
right hand side, without the need for \textit{ad hoc} splitting procedures such
as those customarily employed in most models of this kind and analysed
e.g. in \cite{cullen:2003}, \cite{dubal:2005}, \cite{dubal:2006}.

\section{Numerical results}
\label{tests} 
\indent
A number of numerical experiments have been carried out with  preliminary implementations of the
approach outlined in the previous sections. In particular, since the error bounds presented in \cite{iserles:2001} do not
provide a sharp estimate for the error resulting from the proposed domain decomposition
approach, the goal of these tests is to assess the effective accuracy of the resulting space-time discretization,
as well as to estimate the sensitivity of the results  to the choice of the size of the buffer regions.
Only simple finite difference and finite volume discretizations have been considered in this work,
although it is clear that the results will also depend in general on the chosen spatial discretization.
A further goal of these tests is to understand to which extent the proposed method leads to a reduction
in computational cost with respect to approaches in which the exponential matrix is computed
globally, although,  due to the preliminary nature of the implementation, the estimates
of the CPU times of each method are to be considered  only as rough indications of its computational cost.
In all the tests, the exponential Euler-Rosenbrock methods proposed in \cite{hochbruck:1998} have been employed,
which reduce in the linear case to the exponential Euler method.  In the one dimensional tests,
the $\phi $ matrix was computed by  the scaling and squaring
algorithm and Pad\'e approximation.
For the linear test cases, the $\phi $ matrices associated to each subdomain were computed only
at the first time step. In the nonlinear test cases, the Jacobian of the ODE system
and the corresponding $\phi $ matrices were computed at the first time step and later kept constant
for a number of time steps chosen empirically for each case.

\subsection{One-dimensional, linear tests}
\label{tests_1d_lin}
In a first set of numerical tests, the simple one-dimensional, constant coefficient, linear advection diffusion equation 
\begin{equation}
\label{test1d_1}
  \frac{\partial {c} } {\partial t}=-u\frac{\partial {c} } {\partial x}+ \nu \frac{\partial^2 c } {\partial x^2} 
  \end{equation}
has been considered on a time interval at $[0,T] $  and on an spatial interval $[0,L] $ with periodic boundary conditions,
along with  the one-dimensional Schr\"odinger equation with harmonic potential
\begin{equation}
\label{test1d_2}
 \frac{\partial c} {\partial t}=\frac{i}2 \frac{\partial^2 c } {\partial x^2} -i\frac{\kappa}2 x^2 c
  \end{equation}
  with periodic boundary conditions on an interval $[-L/2,L/2] $ and $\kappa=10.$
In both cases, the differential operators have been approximated by simple centered finite difference formulae
and  a Gaussian initial datum was considered.

\begin{figure}[htbc]
\begin{center}
\includegraphics[height=0.3\textheight,width=0.45\textwidth]{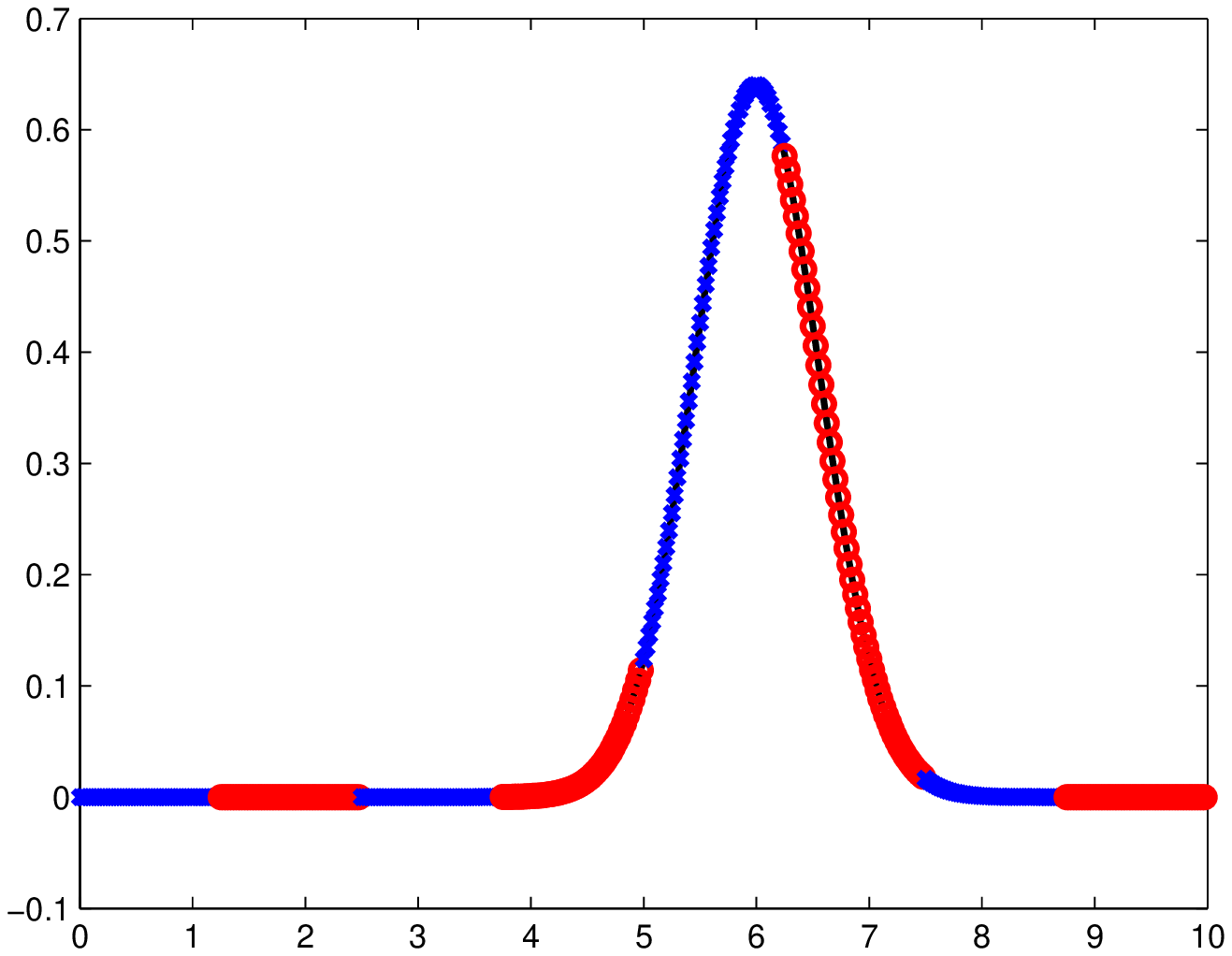}(a)
\includegraphics[height=0.3\textheight,width=0.45\textwidth]{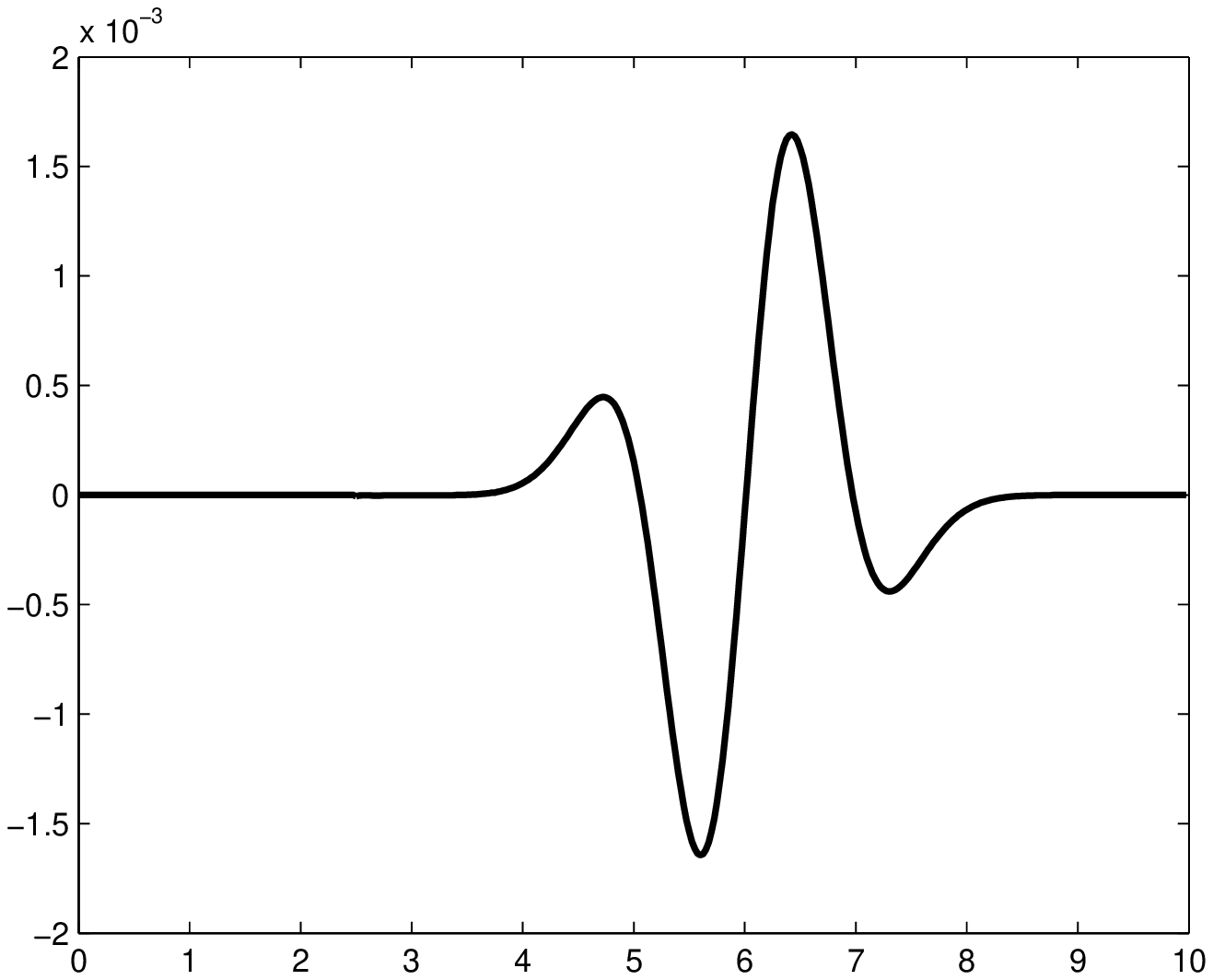}(b)
\end{center}
\caption{Approximation of the advection diffusion equation  at $T=3$ by local exponential method applied over
 8 subdomains with sufficiently large overlap: (a) reference solution by Fourier method and separation of variables (black line) and numerical solution (blue or red symbols depending on the subdomain); (b) absolute error.}
 \label{ddgood}
\end{figure}

\begin{figure}[htbc]
\begin{center}
\includegraphics[height=0.3\textheight,width=0.45\textwidth]{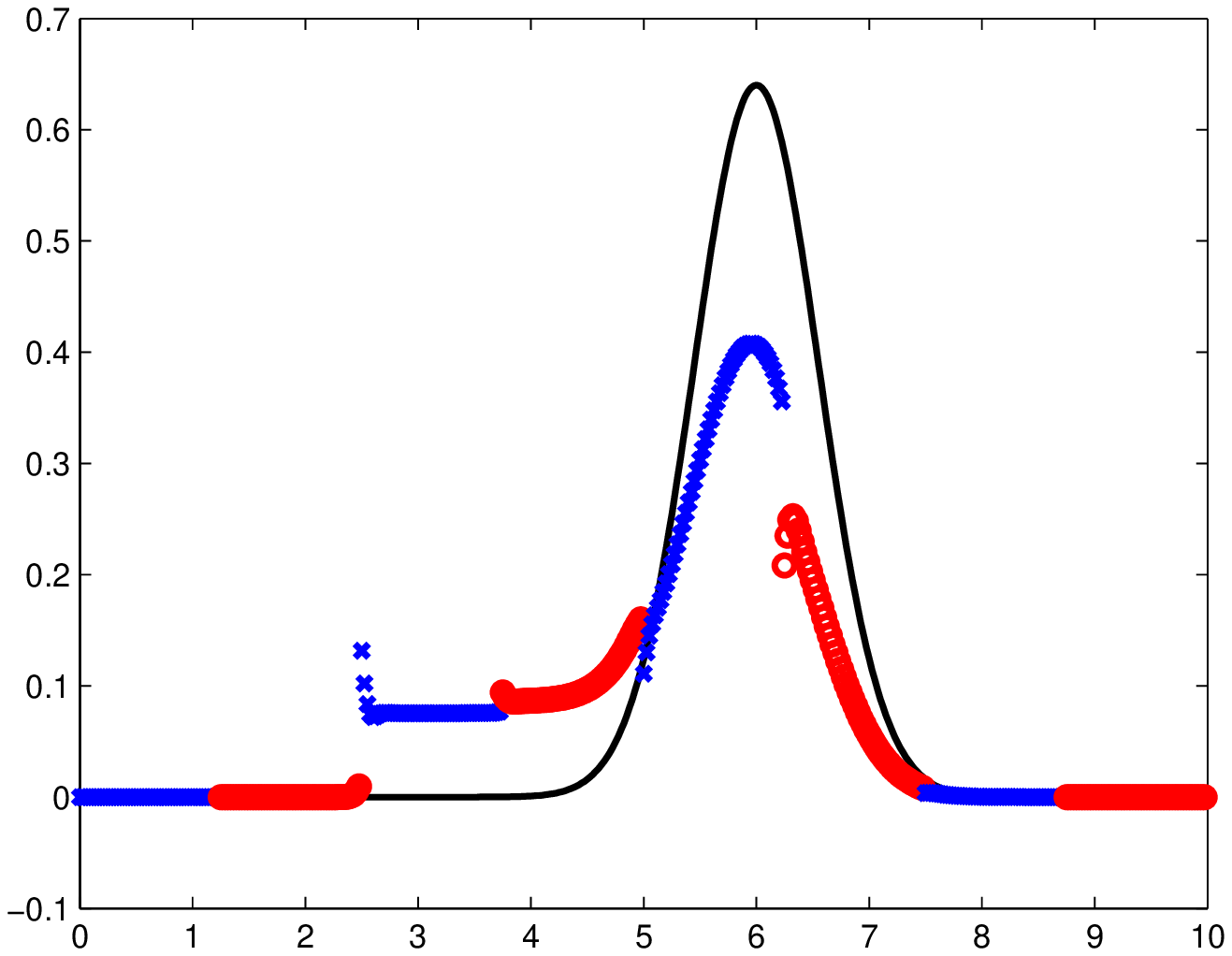}(a)
\includegraphics[height=0.3\textheight,width=0.45\textwidth]{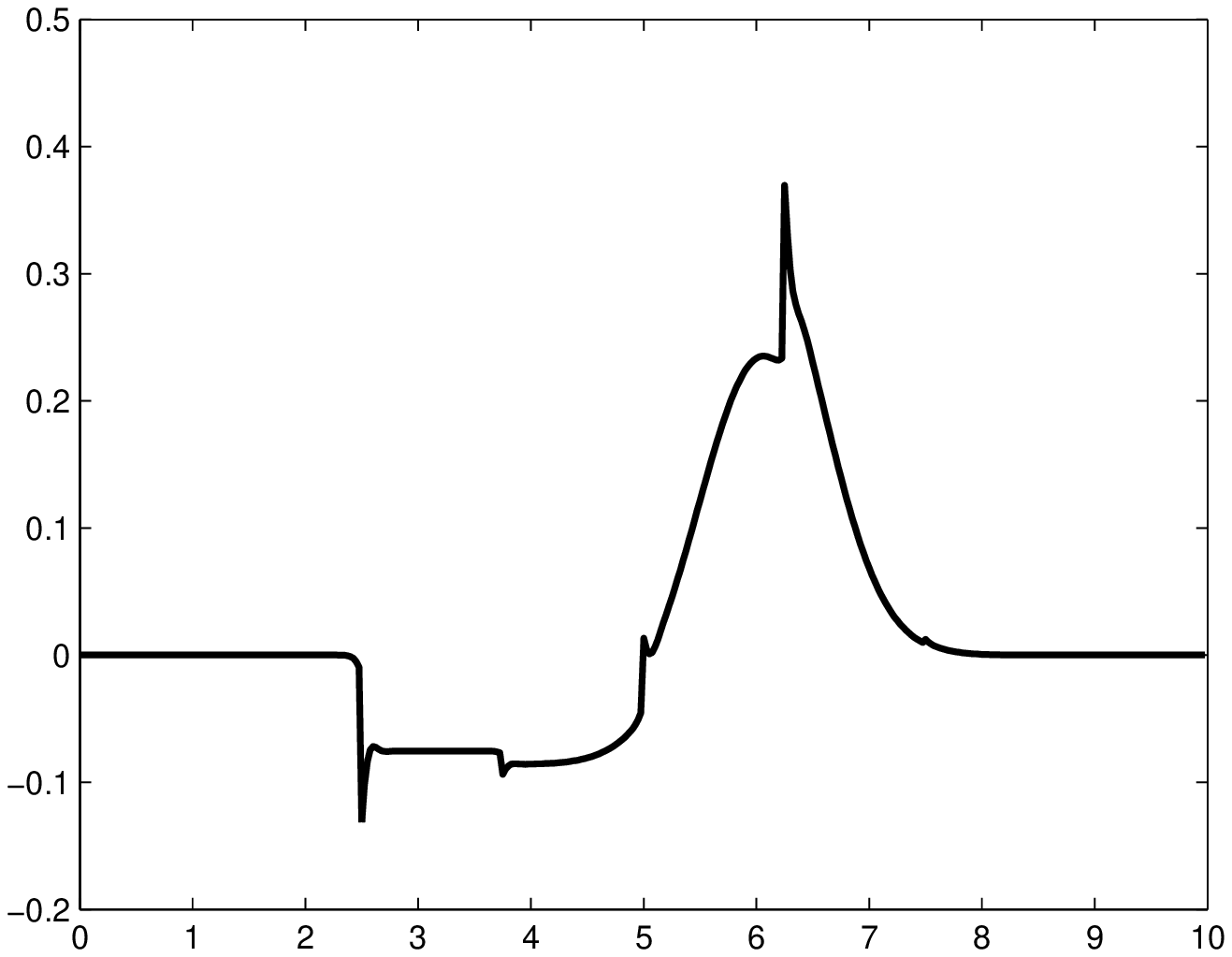}(b)
\end{center}
\caption{Approximation of the advection diffusion equation  at $T=3$ by local exponential method applied over
 8 subdomains with insufficiently large overlap: (a) reference solution by Fourier method and separation of variables (black line) and numerical solution (blue or red symbols depending on the subdomain); (b) absolute error.}
  \label{ddbad}
\end{figure}

The first goal of the tests is to show that the approach outlined in section \ref{dd} does not lead
 in practice to any loss of accuracy, as long
as a sufficiently large overlap is allowed among neighboring subdomains.
 Numerical experiments confirm that this is indeed the case. As an example,
 the numerical solution of \eqref{test1d_1}
  computed at $T=3 $ time units on a domain of size $L=10$
 is displayed in figure \ref{ddgood}. The 
 discretization employs 400 grid points subdivided into 8 identical
 subdomains. The simulation was run with $u, \nu,$ and $\Delta t $ values resulting in the values
 $C\approx 4, $ $ \mu \approx 4.8   $ for the usual
 stability parameters and using buffer regions of size 18 grid points on each side of the subdomains considered.
 It can be seen that the error structure is regular in space and only depends on the spatial derivatives of the solution,
 as expected in the case of a generic time discretizations. As in the case of the standard exponential
 method, the dominant error component is due to the spatial discretization error, as it can be
 seen by comparing the result with a reference solution of the ODE system associated to the
 same  spatial semi-discretization obtained by a high order accurate reference solver with automatic error
 estimation and error tolerance of the order $10^{-9}.$ The errors obtained by the local exponential method
 are, as long as the buffer region is sufficiently large, essentially identical to those of the
 standard exponential method applied by global computation of the exponential matrix.
 On the other hand, the error obviously increases as the buffer region size is reduced, ultimately
 leading to totally erroneous solutions where the subdomain imprinting is clearly visible,
 see as an example  the same solution displayed before  as computed with a buffer region of size
 5 gridpoints in  figure \ref{ddbad}.

In order to assess the potential of the proposed domain decomposition
technique for  reduction of the computational cost associated to exponential methods,
 the same test was repeated progressively increasing the time step and changing the number of subdomains
employed.  The results are displayed in table \ref{lem_advdiff1}, where
$D $ denotes the number of subdomains employed and $B $ the number of grid points  in the buffer regions.
 The case $C=8,B=20,D=20 $ was not run, since the   subdomains would have been of the same 
 size as the buffer regions. In all the tests reported in this table, relative errors
 of approximately $3\times 10^{-3} $ were obtained, which is approximately of the same magnitude as
  that obtained on the same test by an explicit Runge Kutta method of order 3 run at $C=0.5,\mu=0.5.$

\begin{table}[htbc]      
\begin{tabular*}{1\textwidth}{@{\extracolsep{\fill}}||l|c|c|c|c|c|c||} 
\hline  
&$D=1$ & $D=2 $ & $D=4 $ & $D=5 $& $D=10$ & $D=20 $  \\
\hline 
 $C=1,\mu=1, B=8$ & 1.85 & 0.59 & 0.39 & 0.47    & 0.75   & 1.35      \\
\hline
 $C=2,\mu=2, B=12$ & 2.85      & 0.65 & 0.28 &  0.30 & 0.65  & 0.78 \\
 \hline
 $C=4,\mu=4, B=15$ & 2.89      & 1.64  & 0.23  &  0.77 & 0.31  & 0.54\\
 \hline
 $C=8,\mu=8, B=20$ & 3.79      & 1.00 & 0.23 &  0.23 & 0.25  & - \\
\hline
 \end{tabular*} 
\caption{CPU times (in seconds) for LEM runs in the linear advection diffusion test case, as a function of the time step,   the number $D $ of subdomains employed and the number of grid points $B $ in the buffer regions.}
\label{lem_advdiff1}
 \end{table} 
Firstly, it must be observed that the CPU times of the standard exponential methods
are increasing as a function of the stability parameters. This may seem a paradox, since
in this case only one matrix function evaluation is necessary for each run. However, the number of
scaling and squaring steps to be performed in the computation of the $\phi $ matrix increases as $\Delta t $
increases, thus leading  to a larger computational cost in the case of longer time steps.
 The potential of the LEM approach is apparent, since CPU times are reduced by a factor ranging from $5 $ to $20. $
The domain decomposition approach also seems competitive with standard implicit methods,
since for example the CPU time for a Crank Nicolson method run at $C=4,\mu= 4$
performed  without repeating the matrix evaluation is around $0.3$ seconds, but with an error that
is approximately one order of magnitude larger.
On the other hand, in this simple case  the fastest LEM runs still
take approximately 5 times longer than an explicit Runge Kutta method of order 3 run at $C=0.5,\mu=0.5.$

In the case of the Schr\"odinger equation, a numerical solution  of  \eqref{test1d_2} was
 computed at $T=1 $ time units on a domain of size $L=10$ on a mesh
 with 400 grid points
 and a reference solution was computed   discretizing in space by a pseudospectral Fourier
 approach and employing a  high order accurate reference ODE solver with automatic error
 estimation and error tolerance of the order $10^{-9} $ for the time discretization.
 The results are displayed in table \ref{lem_advdiff2}. 
  In all the tests reported in this table, relative errors
 of approximately $6\times 10^{-4} $ were obtained, which is approximately of the same magnitude as
  that obtained on the same test by an explicit Runge Kutta method of order 3 run at $\mu=0.2.$
  Notice that, in this case, also the oscillatory term $i\kappa x^2/2 $ has a major impact on stability
  of explicit methods.
 
 \begin{table}[htbc]      
\begin{tabular*}{1\textwidth}{@{\extracolsep{\fill}}||l|c|c|c|c|c||} 
\hline  
&$D=1$ & $D=2 $ & $D=4 $ & $D=5 $& $D=10$   \\
\hline 
  $\mu=2, B=20$ & 8.21      & 5.74 & 4.56 &  4.45 & 6.32    \\
 \hline
 $\mu=4, B=25$ & 12.83      & 3.92  & 2.65 &  3.08 & 5.04   \\
\hline
 \end{tabular*} 
\caption{CPU times (in seconds) for LEM runs in the Schr\"odinger equation  test  case, as a function of the time step,   the number $D $ of subdomains employed and the number of grid points $B $ in the buffer regions.}
\label{lem_advdiff2}
 \end{table} 
 It can be seen that, in this case, the cost reduction is not as impressive as in the case
 of the advection diffusion problem. However,
 it is to be remarked that, for this test, the standard implementation of the exponential methods
 leads to CPU times that are of the same order of that required by the Crank Nicolson method run
 with the same time step.
 Furthermore, the fastest LEM runs take in this case approximately 5 times less than an explicit Runge Kutta method of order 3 run at 
 $\mu=0.2.$ As a result, the LEM approach appears to be competitive both with standard explicit and implicit methods
 in the case of the one dimensional  Schr\"odinger equation.

\subsection{One dimensional, nonlinear tests}
\label{tests_1d_nonlin}
\begin{figure}[htbc]
\begin{center}
\includegraphics[height=0.4\textheight,width=0.7\textwidth]{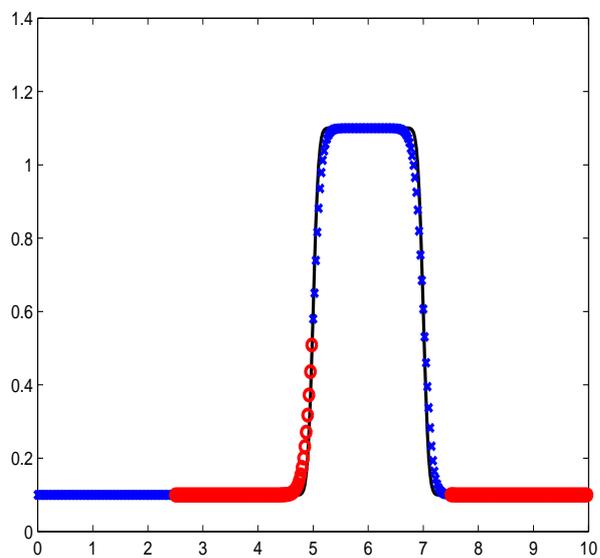} 
\end{center}
\caption{LEM approximation of the solution of the advection diffusion equation discretized  in space by a monotonized finite volume  
method.  The solution is shown at $T=3$ as computed over 
 4 subdomains with sufficiently large overlap:  reference solution (black line) and numerical solution (blue or red symbols depending on the subdomain).}
  \label{ddadvdiffsqwave}
\end{figure}

\begin{figure}[htbc]
\begin{center}
\includegraphics[height=0.4\textheight,width=0.7\textwidth]{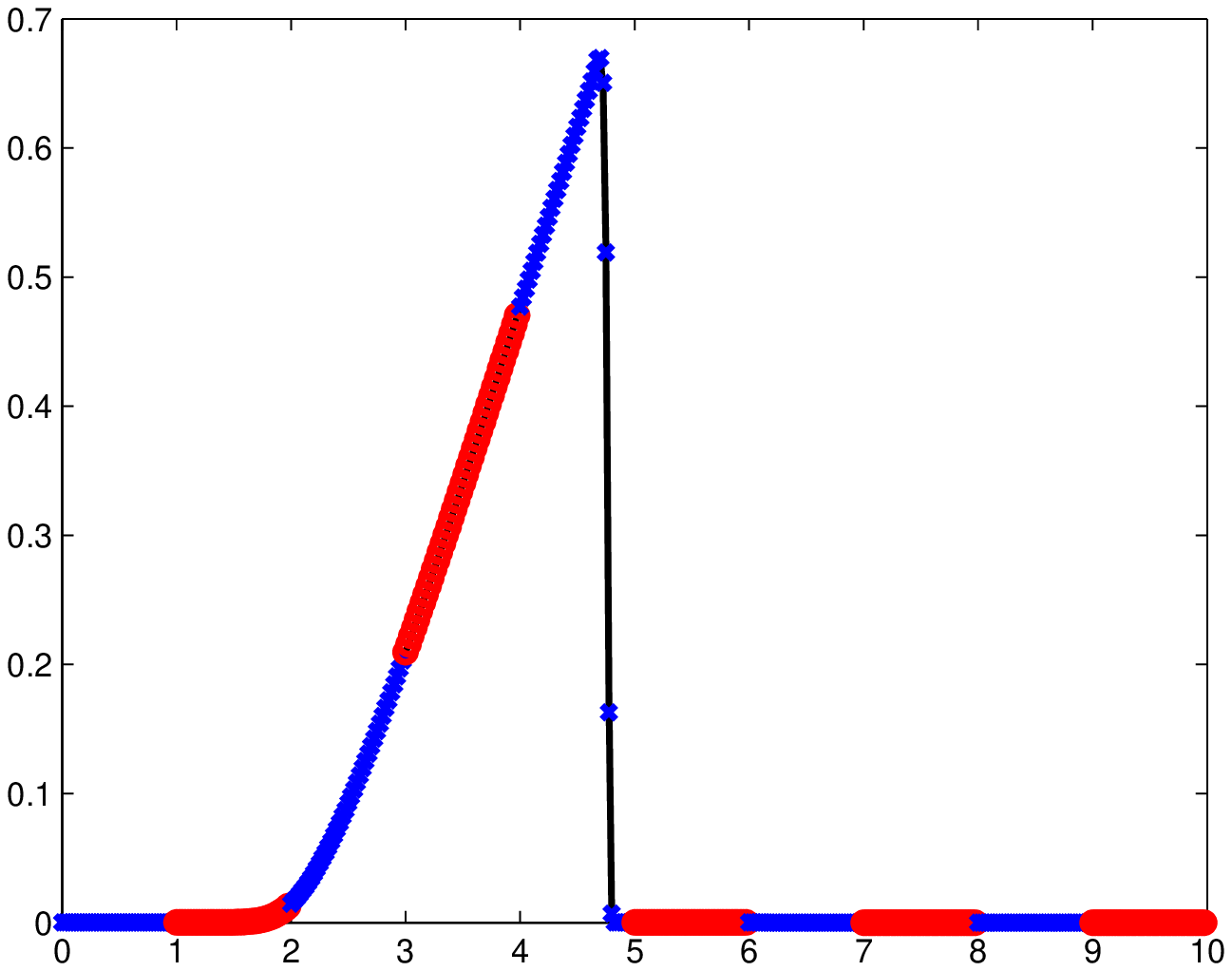} 
\end{center}
\caption{Approximation of the solution of the Burgers equation  at $T=6$ by local exponential method applied over
 10 subdomains with sufficiently large overlap:  reference solution (black line) and numerical solution (blue or red symbols depending on the subdomain).}
 \label{dd_burg_nwave}
\end{figure}

\begin{figure}[htbc]
\begin{center}
\includegraphics[height=0.4\textheight,width=0.7\textwidth]{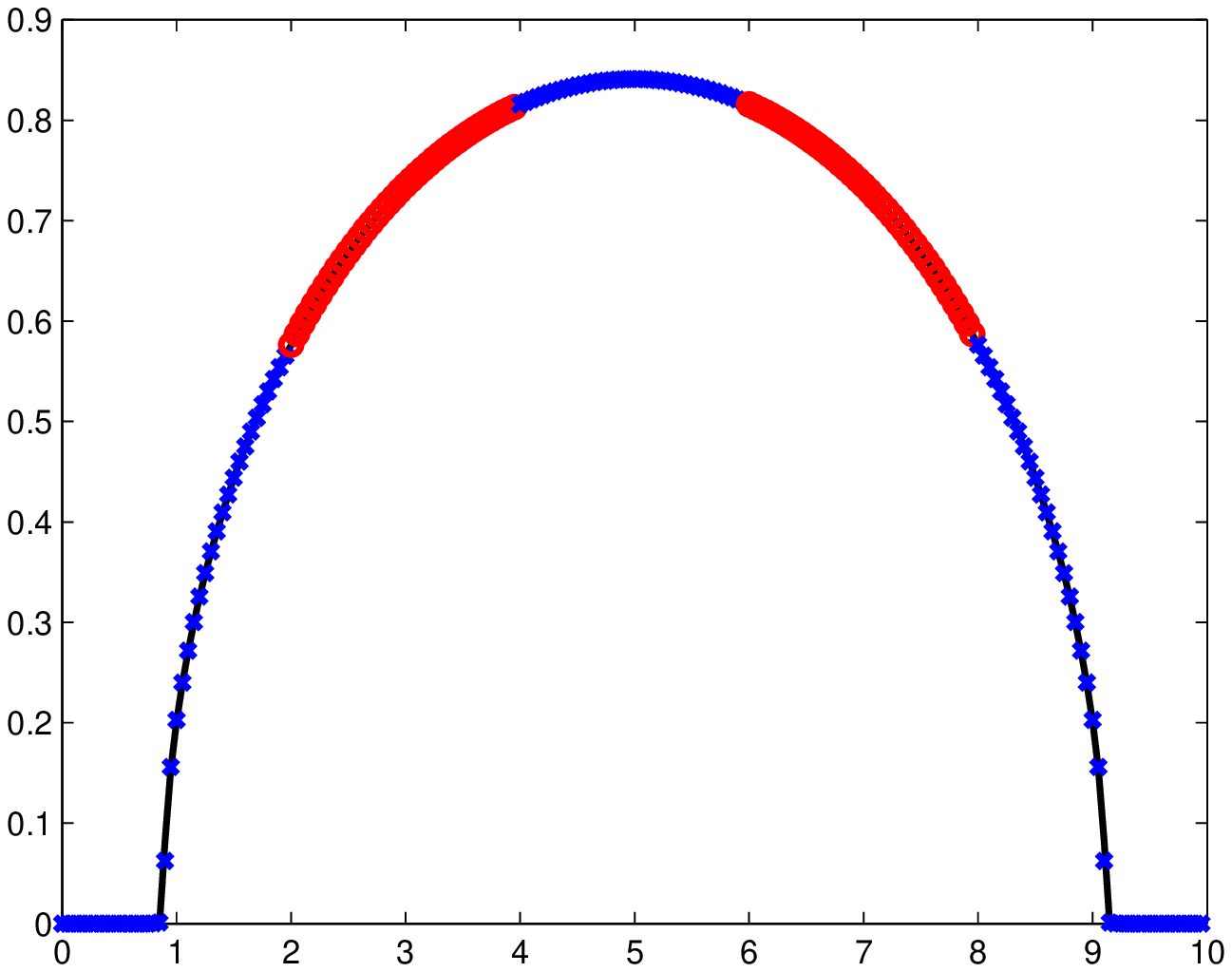} 
\end{center}
\caption{Approximation of the solution of the nonlinear diffusion equation  at $T=1$ by local exponential method applied over 5 subdomains with sufficiently large overlap:  reference solution (black line) and numerical solution (blue or red symbols depending on the subdomain).}
\label{porous_fig}
\end{figure}

Several nonlinear tests have also been performed with the proposed discretization approach.
As a first  nonlinear benchmark,
the time discretization of \eqref{test1d_1} was considered  again, 
with a space discretization given by  a second order, monotonized finite volume approach  employing
a minmod flux limiter. In this case, it is well known (see e.g. \cite{leveque:2002})
that also the space  semi-discretization of a linear advection problem results
in a  nonlinear ODE system. This model problem is relevant for applications
since, in practice, the advection equation is rarely discretized without introducing some analogous
monotonization approach.
Furthermore, an example of  more naturally nonlinear problems, the
 one-dimensional viscous Burgers equation
\begin{equation}
\label{burgers}
\frac{\partial {c} } {\partial t} = -\frac{\partial } {\partial x}\left ( \frac{c^2}2 \right )+\nu \frac{\partial^2 c } {\partial x^2} 
\end{equation}
has been considered, whose nonlinearities are typical of computational
fluid dynamics problems,
along with the nonlinear parabolic equation
\begin{equation}
\label{porous}
  \frac{\partial {c} } {\partial t} =  \frac{\partial^2 c^m } {\partial x^2},
  \end{equation}
for which an exact solution is available (see e.g. \cite{barenblatt:1952}, \cite{bonaventura:2014}), given by
\begin{equation}
u(x,t) = (t+t_0)^{-k} \left(A^2 - \frac{k(m-1)|x|^2}{2m(t+t_0)^{2k}}\right)_+^\frac{1}{m-1}
\end{equation}
where $t_0>0$, $A$ is an arbitrary nonzero constant and
$k = 1/(m+1).$
 A  monotonized second order finite volume approach  was employed
for the spatial discretization also for the Burgers equation,
while equation \ref{porous} was discretized by simple centered finite
differences.  In all cases, an interval of size $L=10$ was considered 
and computational mesh with $400$ control volumes of equal size was employed. Examples  of solutions of these equations obtained by LEM discretization 
employing a third order exponential Rosenbrock method
 are shown in figures \ref{ddadvdiffsqwave}, \ref{dd_burg_nwave}, \ref{porous_fig}, respectively.
 In the case of the advection diffusion equation, a square wave initial datum was considered,
 while in the case of the Burgers equation the initial datum was taken to be Gaussian and for
 equation \ref{porous} the initial condition was chosen so as to recover the analytic
 solution of \cite{barenblatt:1952} with $m=3$, $A=1$.

In a first numerical experiment aimed at checking the performance improvements
in a purely hyperbolic case, the pure advection equation was considered. CPU times
for the LEM discretization with the second order exponential Rosenbrock method
are displayed in table \ref{lem_adv_mon_expros2}, while the corresponding
times for the third order exponential Rosenbrock method
are displayed in table \ref{lem_adv_mon_expros3}. In this test, the final time was taken to be  $T=4$ and
the Jacobian matrices 
used by the Rosenbrock exponential methods and the associated $\phi $ matrices were
recomputed every 5 time steps. Notice that in this test with non smooth
initial datum (and solution), time steps resulting in Courant numbers larger than approximately
1.6 result in violations of   monotonicity for the numerical solution.
In all these tests, the relative $l_{\infty} $ and $l_2$ errors are approximately 0.39 and
0.12, respectively. The error is mostly associated to the spatial discretization error,
as it can be seen comparing to solutions obtained by the same spatial discretization
coupled to a reference solver with small error tolerance. As a comparison,  
 explicit Runge Kutta methods of order 2 and 3 run yield analogous errors when run
 at  Courant numbers between $C=0.2 $ and $C=0.3, $  respectively, with corresponding
 CPU times approximately 0.8 and 1 seconds.
 On the other hand, a Crank-Nicolson time discretization run at Courant numbers
 $C=1 $ and $C=1.6, $   yields relative $l_{\infty} $ errors of 0.45   and 0.5 and CPU times
 of 11.24 and 6.84 seconds, respectively.

\begin{table}[htbc]      
\begin{tabular*}{1\textwidth}{@{\extracolsep{\fill}}||l|c|c|c|c|c|c||} 
\hline  
&$D=1$ & $D=2 $ & $D=4 $ & $D=5 $& $D=8$ & $D=10 $  \\
\hline 
 $C=0.5, B=5$ & 31.7 & 17.24 & 13.00 & 9.72   & 8.81   & 9.11      \\
\hline
 $C=1, B=10$ & 28.77      & 8.96 & 5.59 &  5.46 & 5.47  & 5.48 \\
 \hline
 $C=1.6,  B=15$ & 18.04      & 5.58  & 3.90  &  3.95 & 3.90  & 4.26\\
 \hline
 \end{tabular*} 
\caption{CPU times (in seconds) for LEM runs with second order exponential Rosenbrock method in the advection  test case with monotonized finite volume discretization, as a function of the time step,   the number $D $ of subdomains employed and the number of grid points $B $ in the buffer regions.}
\label{lem_adv_mon_expros2}
 \end{table} 
 
 \begin{table}[htbc]      
\begin{tabular*}{1\textwidth}{@{\extracolsep{\fill}}||l|c|c|c|c|c|c||} 
\hline  
&$D=1$ & $D=2 $ & $D=4 $ & $D=5 $& $D=8$ & $D=10 $  \\
\hline 
 $C=0.5, B=5$ & 32.73 & 15.83 & 10.53 & 9.11   & 9.24   & 9.66      \\
\hline
 $C=1, B=10$ & 15.66      & 8.21 & 5.77 &  5.54 & 5.13  & 6.32 \\
 \hline
 $C=1.6, B=15$ & 9.9      & 6.59  & 4.55  &  4.26 & 4.42  & 5.28\\
 \hline
 \end{tabular*} 
\caption{CPU times (in seconds) for LEM runs with third order exponential Rosenbrock method in the advection  test case with monotonized finite volume discretization, as a function of the time step,   the number $D $ of subdomains employed and the number of grid points $B $ in the buffer regions.}
\label{lem_adv_mon_expros3}
 \end{table} 
 
 In the case of the Burgers equation, a viscosity of $\nu=0.05$.
 A reference solution was computed in this case
 by the same spatial discretization
coupled to a reference solver with small error tolerance. Therefore, in this case only
an estimate of  the time discretization error of the different methods is available.
 CPU times
for the LEM discretization with the second order exponential Rosenbrock method
are displayed in table \ref{lem_burg_expros2}, while the corresponding
times for the third order exponential Rosenbrock method
are displayed in table \ref{lem_burg_expros3}. In this test, the final time was taken to be  $T=5$ and
the Jacobian matrices 
used by the Rosenbrock exponential methods and the associated $\phi $ matrices were
recomputed every 5 time steps.  
 
\begin{table}[htbc]      
\begin{tabular*}{1\textwidth}{@{\extracolsep{\fill}}||l|c|c|c|c|c|c||} 
\hline  
&$D=1$ & $D=2 $ & $D=4 $ & $D=5 $& $D=8$ & $D=10 $  \\
\hline 
 $C=0.4,\mu=0.8, B=8$ & 58.26 & 15.77 & 11.06 &   9.06 & 10.32  & 9.76      \\
\hline
 $C=1,\mu=2, B=15$ & 50.91      &11.46 & 4.70 &  4.49 & 4.53  & 5.03 \\
 \hline
 $C=2,\mu=4, B=20$ & 36.35      & 7.93  & 2.80  &  3.19 & 3.29  & 3.13\\
 \hline
   \end{tabular*} 
\caption{CPU times (in seconds) for LEM runs with second order exponential Rosenbrock method  in the Burgers test case, as a function of the time step,   the number $D $ of subdomains employed and the number of grid points $B $ in the buffer regions.}
\label{lem_burg_expros2}
 \end{table}

\begin{table}[htbc]      
\begin{tabular*}{1\textwidth}{@{\extracolsep{\fill}}||l|c|c|c|c|c|c||} 
\hline  
&$D=1$ & $D=2 $ & $D=4 $ & $D=5 $& $D=8$ & $D=10 $  \\
\hline 
 $C=0.4,\mu=0.8, B=5$ & 25.89 & 13.46 & 10.63 & 10.30    & 10.09   & 10.42      \\
\hline
 $C=1,\mu=2, B=15$ & 30.07      & 9.02 & 4.70 &  4.72 & 4.78  & 5.49 \\
 \hline
 $C=2,\mu=4, B=20$ & 26.77      & 6.57  & 2.87  &  3.24 & 3.05  & 3.13\\
\hline
 \end{tabular*} 
\caption{CPU times (in seconds) for LEM runs with third order exponential Rosenbrock method  in the Burgers test case, as a function of the time step,   the number $D $ of subdomains employed and the number of grid points $B $ in the buffer regions.}
\label{lem_burg_expros3}
 \end{table}

\subsection{Two dimensional tests}
\label{tests_2d}

\begin{figure}[htbc]
\begin{center}
\includegraphics[height=0.4\textheight,width=0.7\textwidth]{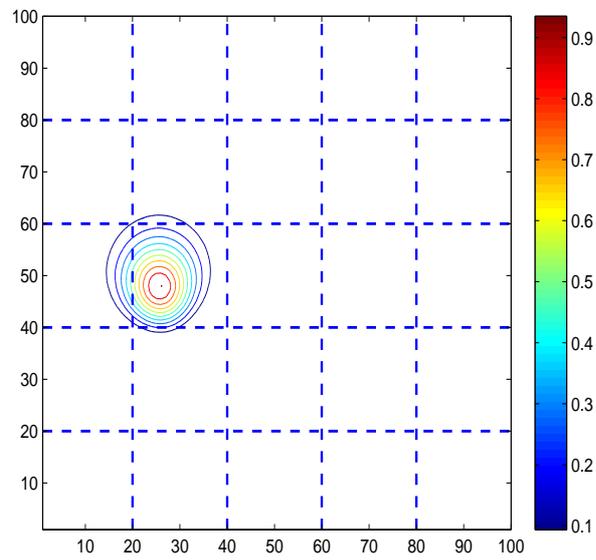} 
\end{center}
\caption{Approximation of the solution of the advection  diffusion equation  in a solid body rotation test. The subdomains employed
are indicated by dashed blue lines.}
\label{advdiff_2d}
\end{figure}

In a preliminary   assessment of the performance of the proposed method
in two dimensions,  an advection diffusion  problem was again
considered, formulated as

\begin{eqnarray}
&&\frac{\partial   c } {\partial t}=- \nabla \cdot ({\bf a}(x)c ) +\nu \Delta   c , \ \ \ \ t \in[0,T]  \ \ \ x \in \Omega \nonumber \\
&&  c(x,0)=  c_0(x) \label{ibv2d}
\end{eqnarray}
where ${\bf a}(x)$ is a divergence free velocity field. As in sections \ref{tests_1d_lin} and \ref{tests_1d_nonlin},
either simple centered finite differences or a monotonized finite volume method were
employed for spatial discretization.
For the approximation of the $\phi $ functions, the Krylov space method of \cite{saad:1992} was employed
for the reference implementation of the standard exponential method. In the case of LEM, the action of the local
 $\phi $ functions was computed, in the present preliminary implementation,  also by the Krylov space method.
 A an example of result in this test case is displayed in figure \ref{advdiff_2d}, where 25 subdomains have been
 employed. Concerning a first quantitative  assessment, the errors of the second order exponential Rosenbrock
 method run at Courant number 7 where analogous to those obtained by an explicit Runge Kutta method of order 4.
 
\begin{figure}[htbc]
\begin{center}
\includegraphics[height=0.4\textheight,width=0.7\textwidth]{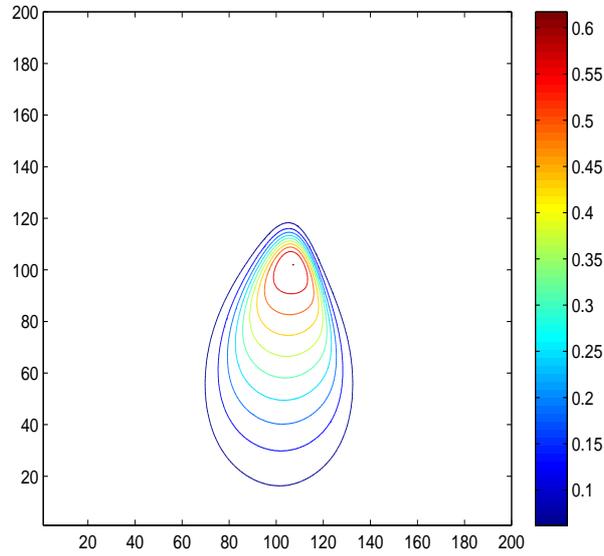} (a)
\includegraphics[height=0.4\textheight,width=0.7\textwidth]{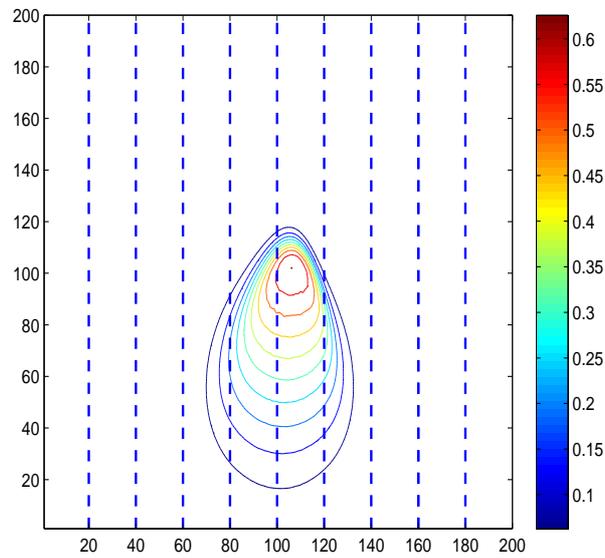}(b) 
\end{center}
\caption{Approximation of the solution of the two-dimensional viscous Burgers equation, computed (a) by a single domain approach
(b) by a multidomain approach. The subdomains employed
are indicated by dashed blue lines.}
\label{burgers_2d}
\end{figure}

 A nonlinear example was also considered, given by a two dimensional extension of a nonlinear Burgers equation.
 This problem was solved on an anisotropic finite volume mesh with vertical spacing much smaller than
 the horizontal one. A result obtained by the second order exponential Rosenbrock
 method, run at Courant number 6 in the vertical direction and Courant number below one in the horizontal
 direction, is shown in figure \ref{burgers_2d}. In this case, a column-wise
 domain decomposition was employed, which allowed to minimize the overlap in the horizontal
 direction.

\section{Conclusions and future work}
\label{conclu}
\indent

An overlapping domain decomposition technique has  been proposed, 
motivated by the results of A.Iserles \cite{iserles:2001},  that allows
to approximate the  global matrices employed in exponential time integration methods for time dependent PDEs
by smaller ones that are related to the spatial subdomains in which the mesh is decomposed.
 The resulting Local Exponential Method (LEM) requires only
the solution of  local problems that can be easily parallelized, thus  increasing the scalability
of the resulting time discretization technique. Furthermore, if the  number
of degrees of freedom associated to each local  
subdomain is small enough, the local exponential matrix can be computed
by simple Pad\'e approximation and can be stored, thus bypassing the
problems that result from having to compute the action  rather than
the exponential matrix itself. The main drawback of the proposed approach is that,  in each of the
local problems, a portion of the local degrees of freedom is only playing an auxiliary role. As a result, their
update is   recomputed multiple times, so that there is an significant overhead with respect to a standard
time discretization,   which is proportional to the size of the buffer regions.
In spite of this overlap, preliminary numerical simulations
show a significant reduction in computational cost with respect to standard
exponential methods.

The results obtained so far appear to justify the further investigation of this approach
in the framework of more complex spatial discretizations and model problems.
Several situations exist in which the proposed approach 
could be useful in spite of its limitation. 
In high order finite element methods, 
high Courant numbers arise easily due to the large number of degrees of freedom per element,
so that a technique that   enables to run at Courant numbers of the order of the polynomial
degree employed with a more local approach should be competitive. Furthermore, in many environmental
applications, such as numerical weather prediction and ocean modelling, strongly anisotropic
meshes are employed, with a vertical resolution that is often  two or three orders of magnitude
smaller than the horizontal one. This results in high vertical Courant numbers, that are often
addressed by directional splitting methods. The present approach would allow to achieve
the same goal by employing a horizontal domain decomposition approach with minimal
overlap among subdomains, such as
almost universally used for parallelization of this kind of models, while at the same time
avoiding \textit{ad hoc} solutions that rely on splitting and providing an efficient 
and robust way to solve  the corresponding fluid dynamics equations.
This would be especially useful for 
models including complex physical parameterizations.
Indeed, LEM would provide  a completely local approach to account for these terms while
maintaining high order accuracy in time without extra computational costs. 
For this reason,  the application of LEM to high order adaptive DG discretizations
will be studied, in order to compare their accuracy and efficiency to that of the semi-implicit, semi-Lagrangian
techniques introduced e.g. in \cite{tumolo:2013}.

\section*{Acknowledgements}
\indent
This work has been partially supported by  the
INDAM - GNCS projects \textit{Metodi ad alta risoluzione per problemi evolutivi fortemente nonlineari} and
 {\it Metodi numerici semi-impliciti e semi-Lagrangiani per sistemi iperbolici di leggi di bilancio},
 as well as by the Office of Naval Research grant  N62909-11-1-4007.
 Several discussions with F.X. Giraldo and M.Restelli on  some of the topics  of this paper are kindly acknowledged.
 I am also grateful to Marco Verani for pointing me to the recent work by Benzi and Simoncini on related problems.
Preliminary results have been presented in the {\it PDE on the sphere 2014} workshop held at NCAR, Boulder, USA,
 in April 2014.

\bibliographystyle{plain}
\bibliography{expint}

\pagebreak

\end{document}